\documentclass[aps,showkeys,10pt]{revtex4}
\usepackage[dvips]{graphicx}
\usepackage{epsfig}
\usepackage{subfigure}
\begin{document}
\title{Integrability of the double pendulum -- the Ramis' question}

\author{Vladimir Salnikov}
\email{vladimir.salnikov@insa-rouen.fr}

\affiliation{Laboratoire de Math\'ematiques de l'INSA de Rouen \\
Avenue de l'Universit\'e
76801 Saint-\'Etienne-du-Rouvray Cedex, 
France }

\begin{abstract}
In this short note we address the problem of integrability of a double pendulum 
in the constant gravity field. We show its non-integrability using the combination 
of algebraic and numerical approaches, namely we compute the non-commuting generators of
the monodromy group along a particular solution obtained numerically. 
\end{abstract}

\maketitle

\section{Introduction}
\label{intro}

A double pendulum is a system that is often modeled as two mass points connected 
between themselves by 
weightless inextensible rods. Mathematically this corresponds to considering the 
motion of two points ${\bf r}_1$ and ${\bf r}_2$ subject to the following constraints:
$$
  \varphi_1 \equiv ({\bf r}_1 - {\bf r}_0)^2 - l_1^2 = 0, \quad
  \varphi_2 \equiv ({\bf r}_2 - {\bf r}_1)^2 - l_2^2 = 0,
$$
the point ${\bf r}_0$ being fixed. 

 \begin{figure}[htp] \centering
    \includegraphics*[height=0.3\linewidth]{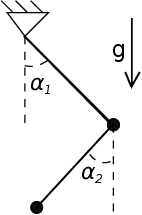}
    \caption{              \label{2pen}
    Double pendulum -- angle parametrization}
 \end{figure}

A convenient approach for the integrability analysis would be resolving this 
constraints, i.e. introducing the generalized coordinates automatically 
taking the constraints into account. As it can be seen from the figure \ref{2pen} 
the system can be parametrized by two angles.
That is for Liouville integrability one needs two independent first integrals. 
The free flat motion of a double pendulum is a classical example of an integrable system: 
the corresponding first integrals are the energy and the  component 
of the angular momentum orthogonal to the plane of the motion. 
But already in the presence of the gravity field the situation is less clear.
It has been shown (mainly using the techniques like the 
Poincar\'e sections) that the system might possess some chaotic behaviour.
The question of Jean-Pierre Ramis sounds roughly as: ``Can one rigorously prove 
non-integrability?''

In this note we discuss meromorphic non-integrability of the double pendulum within the 
framework of the algebraic approach inspired by the results related to 
analysis of the monodromy group (\cite{ziglin}) and the differential Galois group
(\cite{ramis}).
As a mechanical system the double pendulum in the constant gravity 
field $g$ can be described by the following Lagrange function:
\begin{equation} \label{lagr}
  { L} = \dot{\alpha}_1^2  + \dot{\alpha}_1\dot{\alpha}_2 \cos(\alpha_1 - \alpha_2)
  + \frac{1}{2}\dot{\alpha}_2^2 + 2g \cos(\alpha_1) + g \cos(\alpha_2).
\end{equation}
Here and in what follows we assume all the masses as well as all the lengths 
of the constraints to be equal to $1$, this is done to simplify the formula, 
but certainly the analysis can be carried out with arbitrary values of these parameters. 

\section{Description of the method}
The standard procedure of application of the methods of Ziglin or
Morales--Ramis for recovering the obstructions to integrability consists 
of complexifying the equations and 
analyzing the monodromy group or respectively, the differential Galois group, 
of the system of variational equations along an explicitly known particular solution of the 
initial dynamical system. The idea is that for an integrable system 
the obtained group should not be too complicated, i.e. up to some very special cases 
should be (virtually) abelian (see \cite{audin} for details and examples of application).

The key difficulty is that the particular solution should on one hand be 
rather ``nice'' so that one would be able to compute the corresponding monodromy or differential
Galois group, and on the other hand not trivial so that the computed group would contain
a sufficient amount of sources of non-commutativity. 
For example for the system of Euler--Lagrange equations governed by the 
Lagrangian $L$ (eq. \ref{lagr}), that interests us, there are no obvious 
particular solutions except the equilibria. 

To handle this difficulty we have suggested (\cite{ziglin_num}) a method of 
computing the generators of the monodromy group based on a particular solution 
obtained numerically. It consists of basically searching for the paths in the 
complex plane that correspond to finite order branching points of the 
solution of the initial system. Having found several of them we solve the 
system of variational equations in parallel with reconstructing the Riemann surface 
of the particular solution and obtain the corresponding monodromy matrices. 
If among them we see non-commuting ones it permits to conclude meromorphic non-integrability.

\section{Results of the computation}
We are thus solving the system of differential equations 
\begin{equation} \label{syst}
   \dot {\bf x} = {\bf v}({\bf x}), \quad \dot \Xi = A({\bf x}) \Xi, 
\end{equation}
where ``$\dot {\,\,\,\,\,}$'' denotes the derivative with respect to the complex time.
The first equation defines the evolution of 
${\bf x} = (\alpha_1, \alpha_2, \dot {\alpha}_1, \dot {\alpha}_2)^T$ with the right-hand-sides
obtained from the Euler--Lagrange equations governed by $L$ (equation \ref{lagr});
the second one is the matrix system of variational equations with $A(x) = \frac{\partial {\bf v}}{\partial {\bf x}}$.

\begin{figure}[htp] \centering
    \includegraphics*[height=0.4\linewidth]{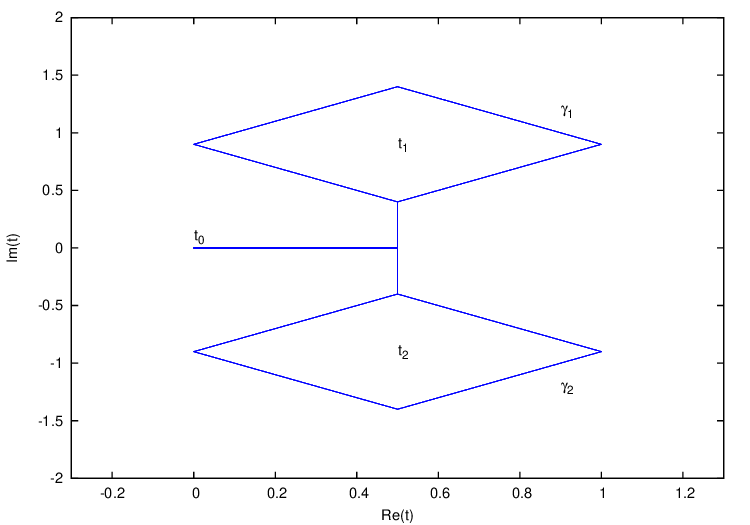}
    \caption{              \label{loops}
    Loops in the complex plane corresponding to the obtained monodromy matrices}
 \end{figure}

Starting from $(\alpha_1, \alpha_2, \dot {\alpha}_1, \dot {\alpha}_2) = (0.1, -0.3, 0.2, 0.4)$
and following the paths $\gamma_1$ and $\gamma_2$ encircling the points $t_1 = 0.5 + 0.9i$ and $t_2 =0.5 - 0.9i$
(three times each of them -- see fig. \ref{loops}) one obtains respectively the monodromy matrices 
$M_1$ and $M_2$.
$$
M_1 = \left(\begin{array}{ccccccc}
20.72 - 17.12i & & 15.79 - 1.34i & &  4.94 + 29.46i & & -14.55 + 10.90i \\ 
-17.67 + 12.78i & & -12.24 - 0.06i & &  -1.93 - 24.87i & & 12.91 - 7.99i \\
12.28 + 6.91i & & 3.55 + 7.78i & &  -13.07 + 7.86i & & -8.18 - 5.42i \\
-11.84 - 12.56i & & -1.31 - 10.39i & &  19.32 - 4.06i & & 8.59 + 9.31i 
      \end{array} \right),
$$
$$
M_2 = \left(\begin{array}{ccccccc}
-18.72 - 17.12i & & -15.79 - 1.34i & &  -4.94413 + 29.46i & & 14.55 + 10.90i \\
17.67 + 12.78i & & 14.24 - 0.06i & &  1.92944 - 24.87i & & -12.91 - 7.99i\\
-12.28 + 6.91i & & -3.55 + 7.78i & &  15.0698 + 7.86i & & 8.18 - 5.42i \\ 
11.84 - 12.56i & & 1.31 - 10.39i & &  -19.319 - 4.06i & & -6.59 + 9.31i  
\end{array} \right).
$$
One easily checks  that these matrices do not commute, that shows meromorphic 
non-integrability. 

It is worth noting that all the computations are done with a controlled 
precision, so the result is certainly more than just a numerical evidence. 
One should also make a remark that the computation passes some obvious compatibility 
tests. For instance if one continues integrating (\ref{syst})
along $\gamma_1$ or $\gamma_2$
for three more loops, one obtains the corresponding matrices squared. 
The eigenvalues are also coherent with their symplectic nature. 
Moreover one can notice, that the obtained matrices can be written 
in the form
$$
  M_1 = Id + A +iB, \quad M_2 = Id - A + iB,
$$
for some real matrices $A$ and $B$, and their commutator 
$[M_1, M_2] = 2i[A, B]$. This observation also simplifies the precision control. 

\section*{Conclusion}
As we have seen in this note the (computer assisted) proof of non-integrability of the 
system describing the motion of a double pendulum is possible. 
It would be interesting to see if the difficulty of the procedure 
depends on parameters of the system and in particular if it is related 
to the value of the energy first integral.

\section*{Acknowledgements}
The author would like to thank the participants of the workshop on Integrability in Dynamical 
Systems and Control -- DISCo--2012 for inspiring discussions related to the subject of this note, 
and in particular Jean-Pierre Ramis for drawing the author's attention to this concrete problem.

\end{document}